\documentclass[a4paper,11pt]{article}

\usepackage{amsmath,color}
\usepackage{amssymb}
\usepackage{graphicx}
\usepackage{ulem}
\usepackage[dvips]{epsfig}
\usepackage[dvips]{graphicx}
\usepackage[dvips]{graphicx}

\definecolor{vertf}{rgb}{0,0.4,0}
\definecolor{violet}{rgb}{0.51,0,0}

\def\build#1_#2^#3{\mathrel{\mathop{\kern 0pt#1}\limits_{#2}^{#3}}}
\newtheorem{prop}{Proposition}

\newtheorem{defi}{Definition}

\newtheorem{theo}{Theorem}

\newtheorem{lem}{Lemma}

\title{The Continuum is Countable: Infinity is Unique}

\author{Laurent Germain\thanks{I thank Herv\'e Boco and also Sylvain Bourjade, Jean Claude Gabillon, Anne Vanhems, Richard for their helpful comments which have improved the quality of this article. I thank also Sandrine, David, Lisa and Th\'eo. Laurent Germain, ESC Toulouse, Universit\'e
de Toulouse, ISAE.
Correspondence to ESC Toulouse, 20
boulevard Lascrosses, BP 7010, 31068 Toulouse cedex 7, France. E-mail:
l.germain@esc-toulouse.fr}}

\begin{document}
\maketitle \baselineskip=0.6cm
\parskip=0.25cm
\thispagestyle{empty}
\newpage
\begin{abstract}
\vspace*{0.2cm}

Since the theory developed by Georg Cantor, mathematicians have
taken a sharp interest in the sizes of infinite sets. We know that
the set of integers is infinitely countable  and that its
cardinality  is $\aleph_0$.  Cantor proved in 1891 with the diagonal
argument that the set of real numbers is uncountable and that there
cannot be any bijection between integers and real numbers.
Cantor states in particular the Continuum Hypothesis. In this paper,
I show that the cardinality of the set of real numbers is the same
as the set of integers. I show also that there is only
one dimension for infinite sets, $\aleph$.     \vspace{0.5cm}

\noindent \textbf{Keywords}: Infinity, Aleph, Set theory.
\vspace{0.5cm}

\vspace*{0.2cm}
\end{abstract}

\setcounter{page}{1}
\newpage

\section{Introduction}
\setcounter{equation}{0} \setcounter{prop}{0} \setcounter{theo}{0}
\setcounter{lem}{0} \setcounter{defi}{0}\setcounter{rem}{0}

Since the theory developed by Georg Cantor, mathematicians have
taken a sharp interest in the sizes of infinite sets. We know that
the set of integers is infinitely countable  and that its
cardinality  is $\aleph_0$.  Cantor proved in 1891 with the diagonal
argument that the set of real numbers is uncountable and that there
cannot be any bijection between integers and real numbers.
 The cardinality of the set of real numbers, the continuum,
 is $c=\aleph_1$. Cantor states, in particular, the Continuum Hypothesis
(CH) by which there is no set the size of which is between the set
of integers and the set of real numbers. The power of the continuum
is equal to $2^{\aleph_0}$ which represents the cardinal of
$P(\textbf{N})$, the set of  all the subsets of $\textbf{N}$. Kurt
Godel in 1939 and Paul Cohen in the 1960s have shown that the CH,
which was mentioned by Hilbert as one of the more acute problems
in mathematics in 1900 in Paris during the International Congress of
Mathematics, was not provable. Paul Cohen showed that the Continuum
Hypothesis is not provable under the Zermelo-Fraenkel set theory
even if the axiom of choice is adopted (ZFC).

In this article, I show that the dimension of the set of integers is
the same as the dimension of the real line. Cantor's theory
mentioned in fact that there were several dimensions for infinity.
This, however, is questionable.  Infinity can be thought as an
absolute concept and there should not exist several dimensions for
the infinite. As a matter of fact, if the set \textbf{N} is an
infinite set it should be of same power as the set \textbf{R}. This
has not been proven so far. We knew by several arguments, and in
particular by the argument of Cantor's diagonal, that the set
\textbf{R} and \textbf{N} do not have the same cardinality and that
there is no bijection between these two sets. In this paper, I show
that there is a bijection between the set of integers and the set of
real numbers. I show that the set \textbf{R} is countable and that
there is only one dimension for infinite sets, $\aleph$. The paper
is organized as follows. In a first section, I show that the set
\textbf{N} can be represented by an infinite tree and that the
cardinal of the set of integers is the same as the cardinal of the
powerset of integers. In a second section, I show that the set of
real numbers can be represented by an infinite tree and that it is
countable. In a third section, I show again that the cardinality of
$P(\textbf{N})$ is the same as the cardinality of \textbf{N}. In a
fourth section, I define infinity. Finally, I make some concluding
remarks.

\section{The Cardinal of the Set of Integers is the Same as the Cardinal of the Powerset of Integers}
\setcounter{equation}{0} \setcounter{prop}{0} \setcounter{theo}{0}
\setcounter{defi}{0} \setcounter{lem}{0}\setcounter{rem}{0}

In this section, I define an infinite tree bijective with \textbf{N} and I show that
the cardinal of $P(\textbf{N})$ is the same as the cardinal of \textbf{N}.

\begin{prop}
\label{pr1} There exists an infinite tree bijective with the set
$\textbf{N}$.
\end{prop}

\noindent \textbf{Proof:} Let us consider an infinite tree starting
with 10 nodes (0,1,2,3,4,5,6,7,8,9). Each of these nodes except the
node 0 gives rise to 10 branches. Each of these 10 branches defines
another 10 nodes. Let us define each of these nodes by the numbers
characterizing the unique path to reach a particular node. The first
node of the tree is the number 0, below 0 there are the numbers
$(1,2,\ldots,9)$. If now we consider the branches that start at node
1 they reach another 10 nodes. These nodes are defined by the
numbers ((10),(11), (12), (13), (14), (15), (16), (17), (18), (19)).
This tree is infinite and any node of the tree gives rise to another
10 nodes. All the numbers of these nodes describe, indeed, the set
$\textbf{N}$. Hence, for any integer there exists a unique path in
the tree and to any node of the tree there corresponds an integer
(see Figure 1 and Figure 2). \textbf{N} is therefore represented by
this infinite tree. The set of the nodes of this infinite tree is
bijective with the set \textbf{N}.

\noindent This proposition shows that there exists an infinite tree
that represents \textbf{N}. In the next proposition, I compute the
cardinality of the powerset of \textbf{N}.

\hspace*{0,6cm}

\begin{prop}
\label{pr2} $\aleph_0$ the cardinal of the set $\textbf{N}$ is equal
to $10^{\aleph_0}$.
\end{prop}

\noindent \textbf{Proof:} Each integer is a node of the tree. 
The set of the nodes in the tree is bijective with 
the set of integers. When N is large, counting the number of nodes is
the same as counting the number of paths (to infinity) in the tree. The
cardinality of the set of the nodes of this infinite tree is
$10^{\aleph_0}$ counted 9 times. Indeed, we see in Figures 1 and 2
that the set $\textbf{N}$ is of cardinality
$(0,1,2,3,4,5,6,7,8,9)^{\textbf{N}}$ counted 9 times. This implies
that the cardinality of the set $\textbf{N}$ is $10^{\aleph_0}$. As
$\aleph_0$ represents the cardinal of $\textbf{N}$, which is by
definition the cardinal of the set of the nodes in the tree
(bijective with $\textbf{N}$), we get $10^{\aleph_0}=\aleph_0$. Note
already that this implies that $2^{\aleph_0}=\aleph_0$.

\noindent This proposition shows that the cardinality of \textbf{N},
$\aleph_0$, is equal to $10^{\aleph_0}$. In the next proposition, I
compute the cardinality of the powerset of \textbf{N}.

\begin{prop}
The cardinal of $P(\textbf{N})$, the powerset of $\textbf{N}$, is
$\aleph_0$.
\end{prop}

\noindent \textbf{Proof:} Note that the tree representing the set
\textbf{N}, which at each node links up with $10$ branches, includes
the infinite subtree with $2$ branches that represents
$P(\textbf{N})$. Indeed, the infinite subtree $(0,1)^{\textbf{N}}$
that represents $P(\textbf{N})$ is included in the tree in bijection
with $\textbf{N}$. This implies that the cardinality of
P(\textbf{N}) is equal to the cardinality of  \textbf{N},
$2^{\aleph_0}=\aleph_0$.

\noindent This proposition shows that the cardinality of the
powerset of \textbf{N} is the same as \textbf{N}. This implies
already that $\aleph_1=\aleph_0$.

\noindent In this section, I have defined an infinite tree which is
bijective with $\textbf{N}$, I have defined the cardinality of
$\textbf{N}$ and shown that the cardinal of the powerset of
$\textbf{N}$ is $\aleph_0$. This section proves already that the power of
the continuum is $\aleph_0$. In the next section, I show the
bijection between integers and real numbers.

\section{Real Numbers}
\setcounter{equation}{0} \setcounter{prop}{0} \setcounter{theo}{0}
\setcounter{lem}{0}\setcounter{defi}{0}\setcounter{rem}{0}

In this section, I define an infinite tree which represents real
numbers, then I define a set of subtrees bijective with $\textbf{N}$
in the interval $(0,1)$. Finally, I compute again the power of the
continuum and the cardinal of the set of real numbers.

\subsection{The Infinite Tree of Real Numbers}

In the next proposition, I define an infinite tree which represents
the set \textbf{R}.

\begin{prop}
\label{pr3} There exists an infinite tree bijective with the set
of all
real numbers defined by their decimal representation.
\end{prop}

\noindent \textbf{Proof:} Let us consider an infinite tree with
$(N+1)$ nodes which give rise to 10 branches, as in Figure 3. Each
of the (N+1) nodes is linked to 10 branches which are again linked
to another 10 branches and so on to infinity. The first $(N+1)$
nodes (in the first column of the tree) represent the set
$\textbf{N}=(0,1,2,\ldots,n,\ldots)$. Each node of this infinite
tree can be defined by the numbers of the branches characterizing
the unique path in the tree to reach a particular node. Thus, the
first node of the tree is the number $0$.  The $10$ other nodes
characterizing the branches that start on the right of the $0$ node
are therefore defined by the couples ((0,1), (0,2), (0,3),
(0,4),(0,5),(0,6),(0,7),(0,8),(0,9),(0,0)) .  For example, the node
(0,1,1) is the node starting at 0 when you take branch 1 and again
the following branch 1 (see Figure 3).

We can now define these nodes by decimal numbers, where the first
figure corresponds to the integer in the first column of the tree
(see Figure 4). For example, the nodes $(0,1), (0,2)$ and $(0,1,1)$
are defined by the decimal numbers $0,1$, $0,2$ and $0,11$. In fact,
we can define any node in this infinite tree by the decimal number
characterizing the path of the tree to reach a node. For example,
the node (3,1,4,1,5,9,2,6…) can be defined by the decimal number
(3,1415926…) which is the number $\pi$. To any integer, or algebraic
or transcendental number, there corresponds a unique path
characterized by a decimal number (see Figure 4). This is true also
for negative numbers that would be drawn on the left of the Figure.
We can notice that, as in Hilbert's Hotel, one can add as many
branches as necessary in the tree to correspond with any sequence of
figures.
 
\noindent This proposition shows that we can establish a bijection
between the infinite tree of Figure 4 and the real numbers. We now
give a lemma that characterize decimal numbers in the tree.

\begin{lem}
If we consider the first $(N+1)$ nodes $(0,1,2,\ldots,n,\ldots)$,
the other $10(N+1)$ branches starting from these nodes define all the 1 decimal real numbers
and the other $(N+1)10^i$ branches define all the i decimal real
numbers, i describing all the set $\textbf{N}$.
\end{lem}

\noindent \textbf{Proof:} For any real number represented by a decimal there exists a
unique path in the tree and vice versa.
We can notice, of course, that in the tree certain real numbers are counted several times.
For example, 0,1 is the same real number as 0,10. The tree represents all the possibilities
to write a real number with decimals.

\noindent Therefore, in this section, I have defined an infinite
tree which represents all real numbers with decimals.

\subsection{Numbers between 0 and 1}

    In this section, I define a set of subtrees that are bijective with $\textbf{N}$ in the interval $(0,1)$.

\noindent     Let us consider all the nodes that depart from the
node 0 and that are between the node 0 and the node 1. Let us define
$t_0$ the particular tree which represents all the paths starting at
the node 0.

\begin{prop} Any real number in the interval $(0,1)$ is defined  by a node of
the infinite tree $t_0$. \end{prop}

\noindent \textbf{Proof:} Any decimal number between 0 and 1 is
characterized by a unique path in the infinite tree (see Figure 6).
For example, let us consider the real number $\pi-3$. There exists a
path in the tree characterizing this particular number.

\noindent We can now establish that to any node of the tree one can
associate an integer.

\begin{prop}
To any node of the tree there corresponds an integer and there
exists a set of subtrees the dimension of which is the same as
$\textbf{N}$.
\end{prop}

\noindent \textbf{Proof:} Let us consider all the nodes that start at the
node 0 and that exclude all the nodes that are linked to the node (00).
One can associate an integer with any node of this truncated
infinite tree. This integer is the number of the node i.e. the
numbers which characterize the unique path to reach this node. Therefore, all
the nodes that start with the number 0 when we exclude all the
branches that depart from the node (0,0) describe the set \textbf{N}.
The application which associates to the nodes of the tree, starting with
the node 0 and excluding the paths starting at the node (0,0), an
 integer  is bijective (see Figures 5, 6, 7 and 8).

\begin{defi} Let $t(0)$ be the set
of all the subtrees $t_{0^i}=t_{(0\ldots0)}$,  0 is counted $i$
times for $i=1$ to $i=N$, where $t_{0^i}$ is the subtree that
represents all the paths that start from the node $(0^i)$ and
exclude all the nodes that are linked to the node $(0^{i+1})$.
\end{defi}

\noindent One can associate an integer with any node of any truncated
infinite tree of $t(0)$. All the nodes that start with the number
$(0^i)$ when we exclude all the branches that depart from the node
$(0^{i+1})$ describe the set \textbf{N}. The application which
associates the nodes of the tree starting with the node $(0^i)$ and
excluding the paths starting at the node $(0^{i+1})$ is bijective
with the set of integers (see Figures 7 and 8). This is true for any
subtree $t_{0^i}$.

\noindent This section defines a set of subtrees bijective with
$\textbf{N}$ in the interval $(0,1)$.

\noindent In the next section, I compute again the power of the continuum.

\subsection{The Power of the Continuum}

In this section, I compute again the power of the continuum showing
a bijection between integers and real numbers.

\begin{theo} The dimension of the interval (0,1) of real numbers is
$\aleph_0$ the dimension of \textbf{N}. \end{theo}

\noindent We know that all the subtrees $t_{0^i}$ for all i and the
set \textbf{N} have the same cardinality. There are $\aleph_0$
subtrees in the set $t(0)$ as there are $\aleph_0$ nodes starting
with $(0^i)$, $i$ describing the set $\textbf{N}$
(see Figures 7 and 8). As the countable union of countable sets is
countable, we get that the cardinality of $t(0)$ is $\aleph_0$.
Hence, there is a bijection between the set of integers and the
interval (0,1). Therefore, the cardinality of the set of real
numbers is indeed $\aleph_0$. We can therefore state the following
theorem.

\begin{theo}

 The power of the continuum is $\aleph=\aleph_0$.
\end{theo}

\noindent \textbf{Proof:} By the same reasoning as before, we find
that there are $\aleph_0$ real numbers in the infinite tree when
considering the union of all the intervals $(0,1),
(1,2),\ldots(n-1,n)$. The power of the continuum is indeed $\aleph=\aleph_0$.

\noindent This section proves again that the continuum is countable.
In the next section, I compute the cardinality of \textbf{R}.

\subsection{Cardinal of Real Numbers}

In this section, I compute again the cardinal of the set of real numbers.

\begin{prop}

The cardinality of \textbf{R}~is $\aleph 10^{\aleph}=\aleph$.
\end{prop}

\noindent \textbf{Proof:}  The cardinality of
the set of the nodes of the infinite tree which represents real 
numbers is ${\aleph}10^{\aleph}$,
that is to say $(0,1,2,3,4,5,6,7,8,9)^\textbf{N}$ counted $\aleph$
times (see Figure 4). As $10{\aleph}=\aleph$ (Proposition 2.2), the
cardinality of \textbf{R} is equal to $\aleph \times \aleph$ which
is equal to $\aleph$.

\noindent This section defines another way to compute the
cardinality of the set of real numbers. In the next section, I
compute again the cardinality of $P(\textbf{N})$.

\section{ Cardinality of the set P(\textbf{N})}

In this section, I compute again the cardinality of the powerset of $\textbf{N}$.

\begin{theo}
$\aleph_1=\aleph_0=\aleph$.
\end{theo}

\noindent \textbf{Proof:} $P(\textbf{N})$ is the powerset of
$\textbf{N}$, the set of all subsets of \textbf{N}. The cardinality of
this set of all subsets is $2^{\aleph_0}$ equal to $\aleph_0=\aleph$, the
power of the continuum. Hence the cardinality of $P(\textbf{N})$ is
$\aleph$.

\noindent This section shows again that continuum is countable. In
the next section, I define infinity.

\section{Unique Infinity}

This section defines the dimension of infinity.

\begin{theo} There is a unique dimension for the infinity which is $\aleph$.
\end{theo}
\noindent \textbf{Proof} We know that the dimension of
$P(\textbf{N})$ is $2^{\aleph_0}$. Therefore, since we know that
$2^{\aleph_0}=\aleph_0$, by induction showing
that $\aleph_{\alpha+1}=\aleph_{\alpha}=\aleph_0=\aleph$ is straightforward.

\noindent As a consequence, there is only one dimension for
infinity.

\noindent This section shows that there is one unique dimension for
infinity.

\section{Conclusion}

In this article, I prove that the cardinality of infinite sets is always
$\aleph$. There is a unique dimension for infinity. I
also prove that infinity is always countable. The
consequence of this result is that the continuum is a countable set.
This result has several consequences in Mathematics, Probability and
Statistics. It modifies not only our vision of the world, but
also that of modeling in Physics, Economics,
Biology and Computer Science, among other fields. Moreover,
it opens  the door to new
concepts in Philosophy.

\begin{center}
\begin{figure}[h]
\epsfxsize=12cm
$$\epsfbox{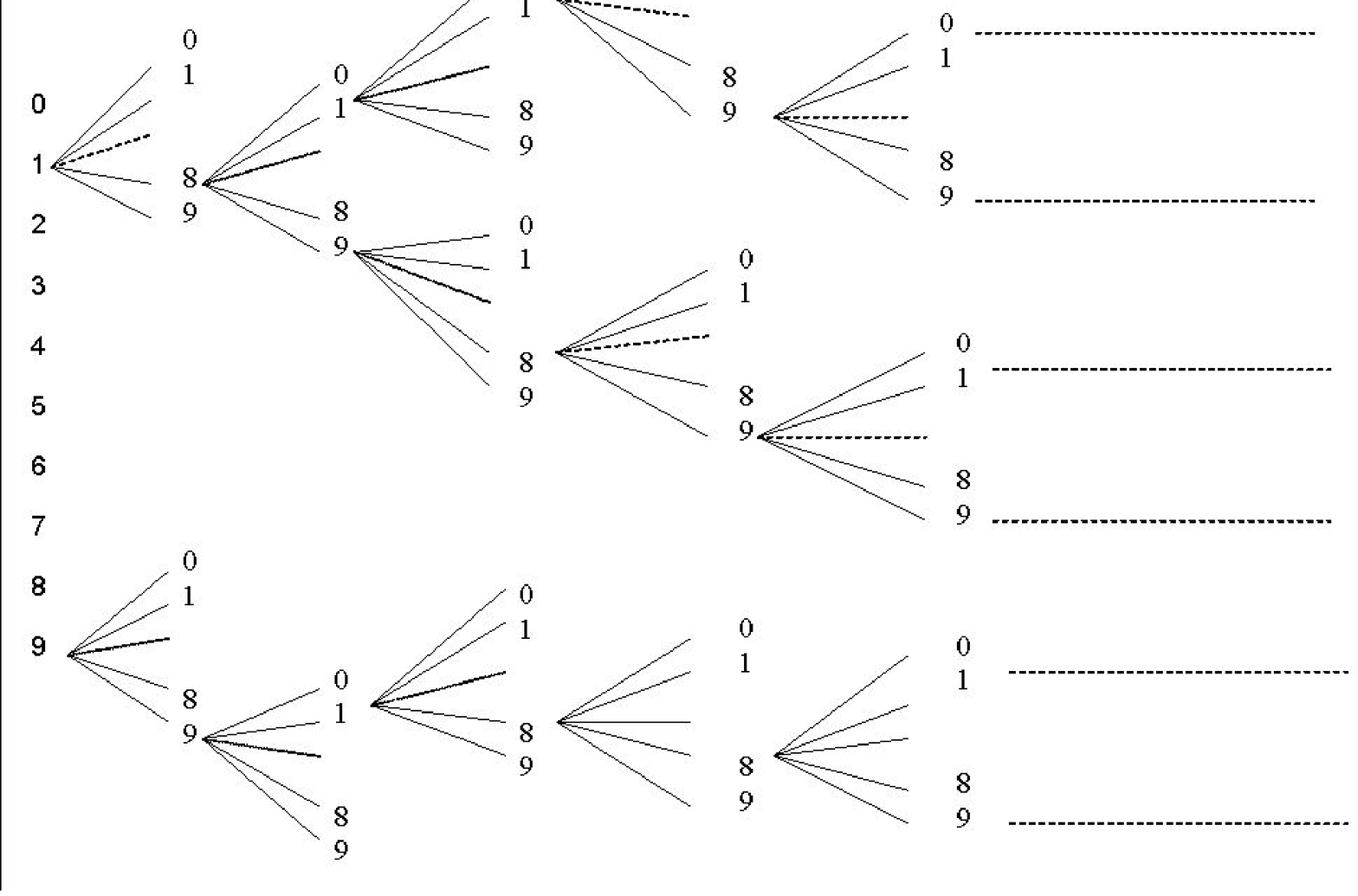}$$\caption{\label{fig7}This figure shows
the infinite tree where each node gives rise to 10 branches.}
\end{figure}
\end{center}

\begin{center}
\begin{figure}[h]
\epsfxsize=12cm $$\epsfbox{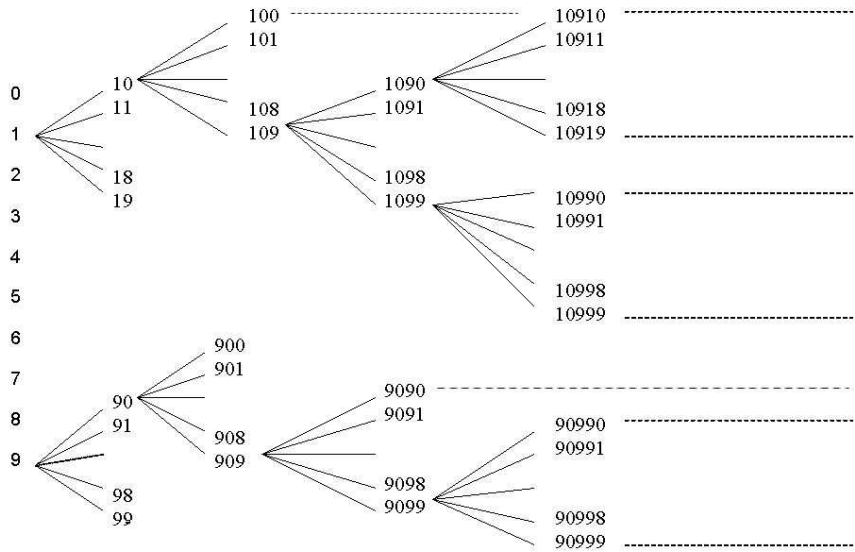}$$ \caption{\label{fig7}This
figures shows the infinite tree representing the set \textbf{N}.}

\end{figure}
\end{center}

\begin{center}
\begin{figure}[h]
\epsfxsize=12cm $$\epsfbox{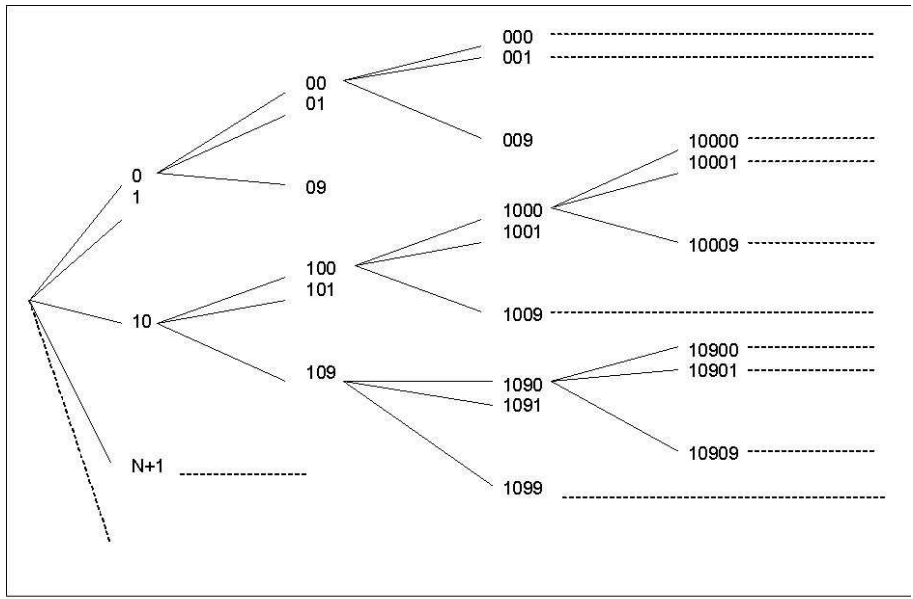}$$ \caption{\label{fig7}This
figure shows the infinite tree where each node is defined by the
unique path to reach a node.}
\end{figure}
\end{center}

\begin{center}
\begin{figure}[h]
\epsfxsize=12cm $$\epsfbox{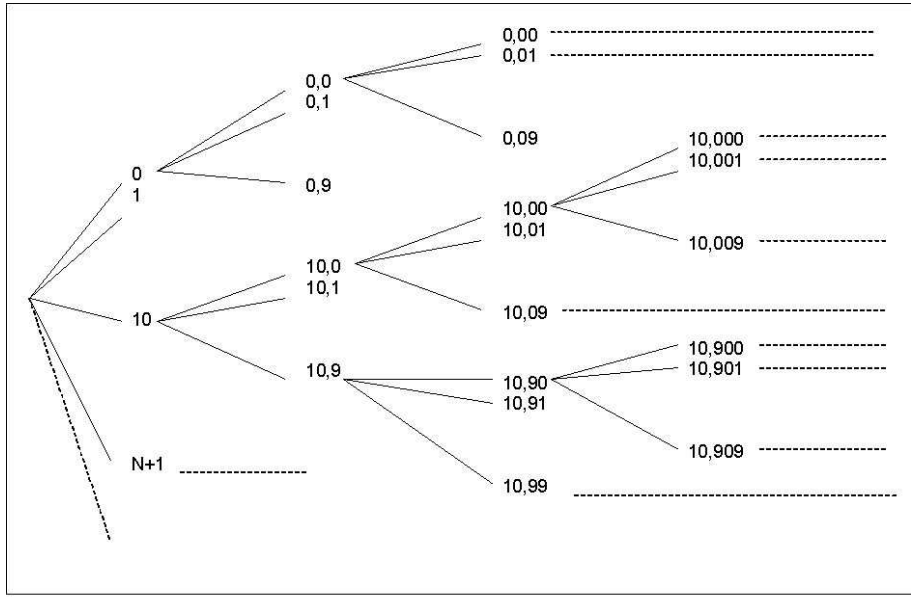}$$ \caption{\label{fig7}This
figure shows the infinite tree where  all real
numbers are represented.}
\end{figure}
\end{center}

\begin{center}
\begin{figure}[h]
\epsfxsize=12cm $$\epsfbox{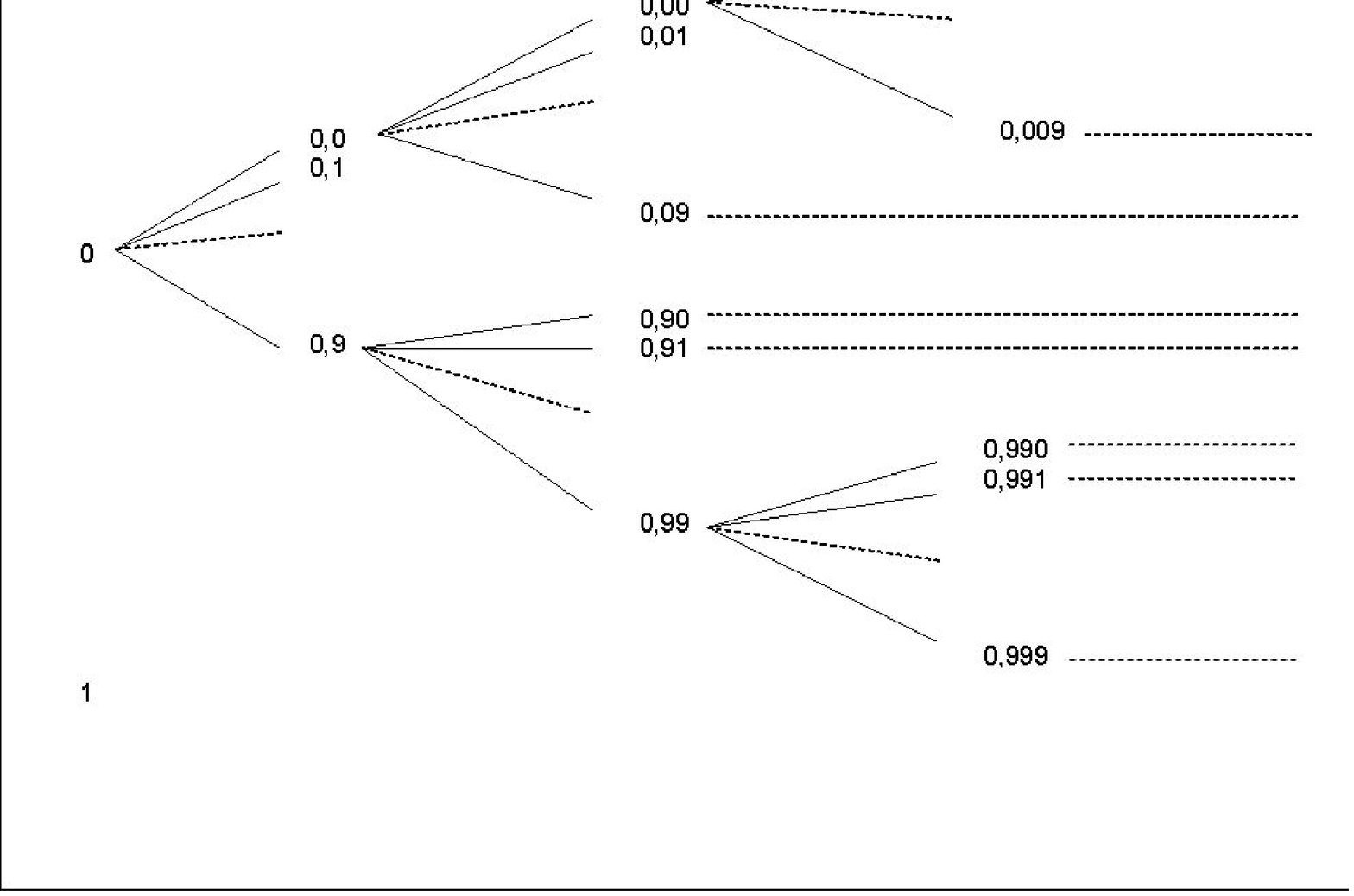}$$ \caption{\label{fig7}This
figure shows the infinite tree representing all the real numbers in
the interval (0,1).}
\end{figure}
\end{center}

\begin{center}
\begin{figure}[h]
\epsfxsize=12cm $$\epsfbox{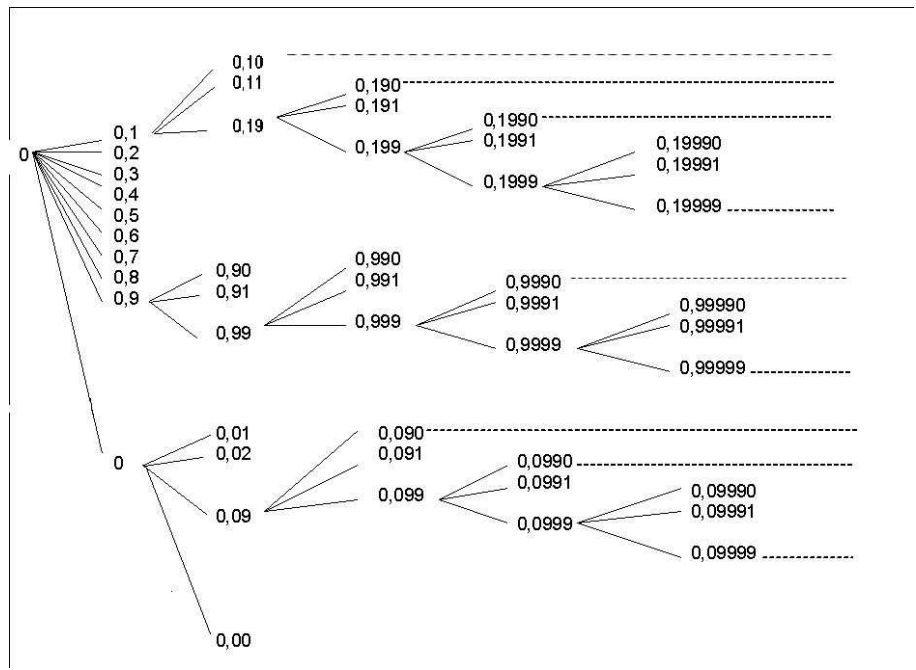}$$ \caption{\label{fig7}This
figure shows the infinite tree of the interval (0,1) rearranged with
nodes starting with zeros in the downside.}
\end{figure}
\end{center}

\begin{center}
\begin{figure}[h]
\epsfxsize=12cm $$\epsfbox{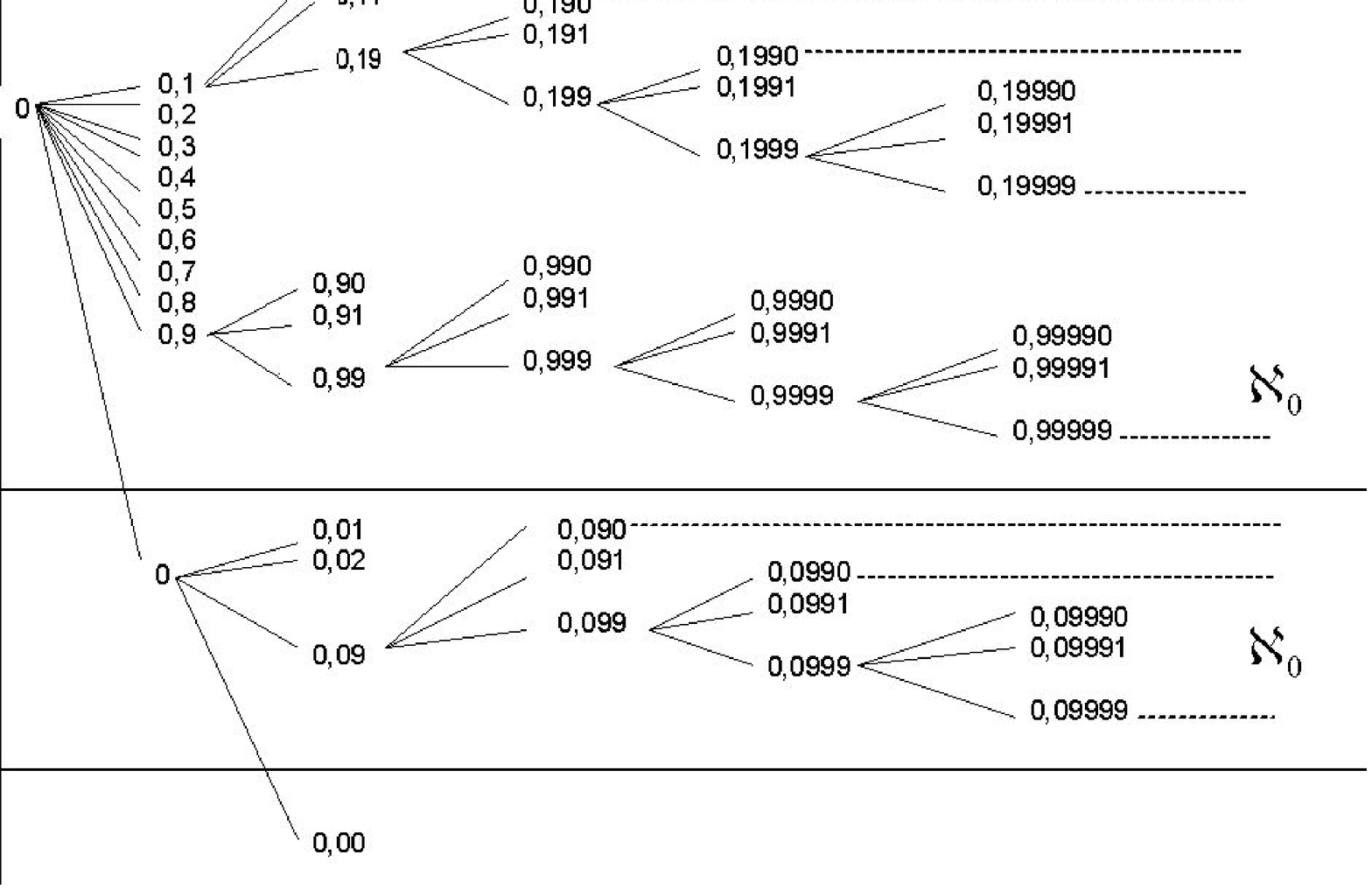}$$ \caption{\label{fig7}This
figure shows the truncated infinite trees   of
dimension $\aleph_0$.}
\end{figure}
\end{center}

\newpage

\begin{center}
\begin{figure}[h]
\epsfxsize=12cm $$\epsfbox{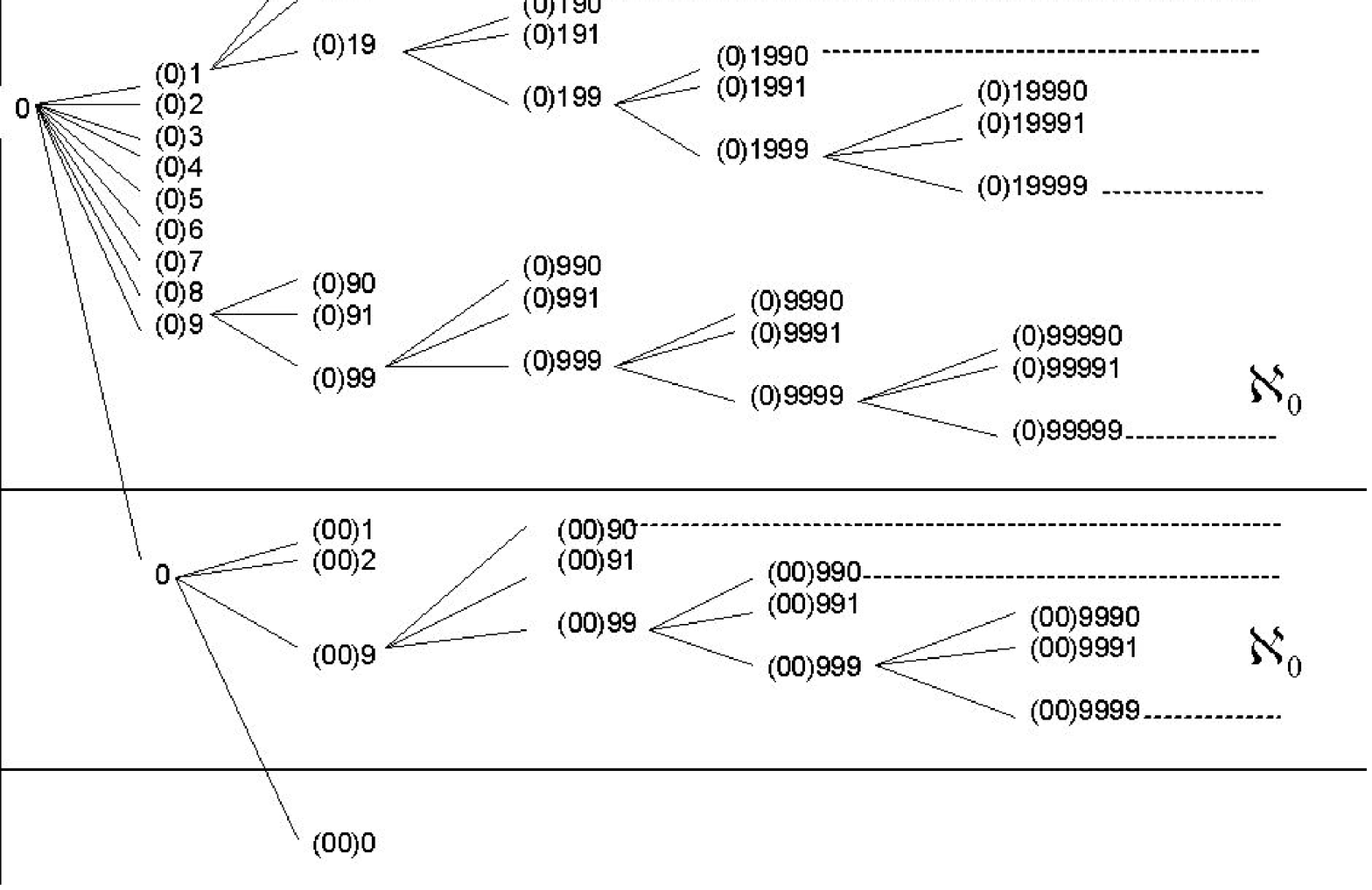}$$ \caption{\label{fig7}This
figure shows the bijection between the truncated trees and
\textbf{N}.}
\end{figure}
\end{center}


\begin{center}
{\Large \textbf{References} }
\end{center}


\begin{enumerate}
\setlength{\parskip}{-1.5mm}

\item Cantor, G.: Uber eine elementare Frage zur Mannigfaltigkeitslehre, Jahresbericht der Deutschen Mathematiker-Vereinigung (1891).

\item Cohen, Paul J.: The Independence of the
Continuum Hypothesis. Proceedings of the National Academy of
Sciences of the United States of America (1963).

\item G¨odel, K.: The Consistency of the Continuum-Hypothesis.
Princeton University Press (1940).

\item  G¨odel,K.: What is Cantor\'~s Continuum Problem?. Amer. Math.
Monthly 54, 515$-$525 (1947).

\end{enumerate}
\end{document}